\documentclass[a4paper,10pt]{article}
\usepackage{amsfonts}
\usepackage{amstext}
\usepackage{amsthm}
\usepackage{amsmath}
\usepackage{amssymb}
\usepackage{latexsym}
\usepackage[latin1]{inputenc}
\newcommand{\bl}{\hfill\rule{2mm}{2mm}}
\newcommand{\R}{\Bbb{R}}

\newtheorem{teor}{Theorem}[section]
\newtheorem{propo}{Proposition}[section]

\newtheorem{lema}{Lemma}[section]

\newcommand{\n}{\noindent}

\begin{document}

\title{Sharp $L^p$-entropy inequalities on manifolds
 }

\author{\textbf{Jurandir Ceccon \footnote{\textit{E-mail addresses}:
ceccon@ufpr.br (J. Ceccon)}}\\ {\small\it Departamento
de Matem\'{a}tica, Universidade Federal do Paran\'{a},}\\ {\small\it Caixa Postal 19081, 81531-980, Curitiba, PR, Brazil}\\
\textbf{Marcos Montenegro \footnote{\textit{E-mail addresses}:
montene@mat.ufmg.br (M. Montenegro)}}\\ {\small\it Departamento de Matem\'{a}tica,
Universidade Federal de Minas Gerais,}\\ {\small\it Caixa Postal 702, 30123-970, Belo Horizonte, MG, Brazil}} \maketitle

\markboth{abstract}{abstract}
\addcontentsline{toc}{chapter}{abstract}

%\hrule \vspace{0,2cm}

\begin{center}
\small{
{\bf Abstract}}
\end{center}

\small{In 2003, Del Pino and Dolbeault \cite{DPDo} and Gentil
\cite{G} investigated, independently, best constants and extremals
associated to Euclidean $L^p$-entropy inequalities for $p > 1$. In this
work, we present some contributions in the Riemannian
context. Namely, let $(M,g)$ be a closed Riemannian manifold of
dimension $n \geq 3$. For $1 < p \leq 2$, we establish the validity of the sharp
Riemannian $L^p$-entropy inequality

\[
\int_M |u|^p \log(|u|^p) dv_g \leq \frac{n}{p} \log \left( {\cal A}_{opt} \int_M |\nabla_g u|^p dv_g + {\cal B} \right)
\]

\n on all functions $u \in H^{1,p}(M)$ such that
$||u||_{L^p(M)} = 1$ for some constant ${\cal B}$. Moreover, we prove that the first best constant ${\cal A}_{opt}$ is equal to the corresponding
Euclidean one. Our approach is inspired on the Bakry, Coulhon,
Ledoux and Sallof-Coste's idea \cite{Ba} of getting Euclidean
entropy inequalities as a limit case of suitable subcritical interpolation inequalities. It is conjectured that the inequality sometimes fails for $p > 2$.}

\vspace{0,5cm}
%\hrule\vspace{0.2cm}

\begin{center}
\section{Introduction}
\end{center}

Logarithmic Sobolev inequalities are a powerful tool in Real
Analysis, Complex Analysis, Geometric Analysis, Convex Geometry
and Probability. The pioneer work by L. Gross \cite{Gr} puts
forward the equivalence between a class of Euclidean logarithmic Sobolev
inequalities and hypercontractivity of the associated heat
semigroup. Particularly, his logarithmic Sobolev inequality with respect to the Gaussian measure plays an important role in Ricci flow theory ({\it e.g.} \cite{Pe}), optimal transport theory ({\it e.g.} \cite{Tal}), probability theory ({\it e.g.} \cite{Le1}), among other applications. Later, Weissler \cite{W} introduced a log-Sobolev type inequality (also known as Euclidean $L^2$-entropy inequality) equivalent to the Gross's inequality with Gaussian measure. The Euclidean $L^2$-entropy inequality and its variants have been used in the study of optimal estimates for solutions of certain nonlinear diffusion equations, see for instance \cite{BGL}, \cite{BGL1}, \cite{Bro1}, \cite{DDG}, \cite{G0}, \cite{G} and references therein.

The Euclidean $L^p$-entropy inequality for $p \geq 1$ states that, for any function $u \in W^{1,p}(\R^n)$ with $\int_{\R^n} |u|^p dx = 1$,

\begin{equation}\label{dee}
Ent_{dx}(|u|^p):= \int_{\R^n} |u|^p \log(|u|^p)\; dx \leq \frac{n}{p} \log \left({\cal A}_0(p) \int_{\R^n} |\nabla u|^p\; dx\right)\; ,
\end{equation}

\n where $n \geq 2$, $p \geq 1$ and ${\cal A}_0(p)$ is the best possible constant in this inequality.

As mentioned above, the Euclidean $L^2$-entropy inequality was established by Weissler in \cite{W}. Thereafter, Carlen \cite{Carlen} showed that its extremal functions are precisely dilations and translations of the Gaussian function

\[
u_0(x) = \pi^{-\frac n2} e^{-|x|^2}\; .
\]

\n For $p = 1$, Ledoux \cite{Le} proved the inequality (\ref{dee}) and Beckner \cite{Be} classified its extremal
functions as normalized characteristic functions of balls. In \cite{Be}, Beckner also proved that (\ref{dee}) is valid for any $1 < p <
n$ and Del Pino and Dolbeault \cite{DPDo} characterized its extremal functions as dilations and translations of the function

\begin{equation} \label{ext}
u_0(x) = \pi^{-\frac n2} \frac{\Gamma(\frac{n}{2} +
1)}{\Gamma(\frac{n(p - 1)}{p} + 1)} e^{-|x|^{\frac{p}{p-1}}}\; .
\end{equation}

\n Finally, Gentil \cite{G} established the validity of (\ref{dee}) and that $u_0$ is an extremal function for any $p > 1$. Thanks to an uniqueness argument due to Del Pino and Dolbeault \cite{DPDo}, modulo dilations and translations, the classification also extends for any $p > 1$.

As a byproduct, they derived, for any $p > 1$,

\[
{\cal A}_0(p) = \frac{p}{n} \left( \frac{p - 1}{e} \right)^{p - 1}
\pi^{-\frac{p}{2}} \left( \frac{\Gamma(\frac{n}{2} +
1)}{\Gamma(\frac{n(p - 1)}{p} + 1)} \right)^{\frac{p}{n}}\; .
\]

In order to introduce sharp $L^p$-entropy inequalities within the Riemannian environment for $1 \leq p < n$, we first deduce an intermediate entropy inequality.

Let $(M,g)$ be a smooth closed Riemannian manifold of dimension $n \geq 2$ and $1 \leq p < n$. For any $p \leq s \leq p^*:=\frac{np}{n-p}$, H\"{o}lder's inequality gives

\[
||u||_{L^s(M)} \leq ||u||_{L^p(M)}^\alpha ||u||_{L^{p^*}(M)}^{1 -
\alpha}\; ,
\]

\n for all function $u \in L^{p^*}(M)$, where $\alpha = \frac{np - ns +
ps}{ps}$. Taking logarithm of both sides, one gets

\[
\log \left( \frac{||u||_{L^s(M)}}{||u||_{L^p(M)}} \right) +
(\alpha - 1) \log \left( \frac{||u||_{L^p(M)}}{||u||_{L^{p^*}(M)}}
\right) \leq 0\; .
\]

\n Since this inequality trivializes to an equality when $s = p$,
we may differentiate it with respect to $s$ at $s = p$. Then, a simple computation provides

\[
\int_M |u|^p \log |u|^p \; dv_g \leq \frac{n}{p} \log \left(
\int_M |u|^{p^*} \; dv_g\right)^{\frac{p}{p^*}}
\]

\n for all function $u \in L^{p^*}(M)$ with $||u||_{L^p(M)} = 1$.

Using the above inequality and the Sobolev embedding theorem for compact manifolds ({\it e.g.} \cite{Au4}), one gets constants ${\cal A}, {\cal B} \in \R$ such
that, for any $u \in H^{1,p}(M)$ with $\int_M |u|^p dv_g = 1$,

\begin{gather}\label{AB}
Ent_{dv_g}(|u|^p):= \int_{M} |u|^p \log |u|^p\; dv_g \leq \frac{n}{p} \log \left( {\cal A} \int_M |\nabla_g u|^p\; dv_g + {\cal B} \right)\; , \tag{$L({\cal A},{\cal B})$}
\end{gather}

\n where $dv_g$ denotes the Riemannian volume element, $\nabla_g$
is the gradient operator of $g$ and $H^{1,p}(M)$ is the Sobolev
space defined as the completion of $C^{\infty}(M)$ under the norm

\[
||u||_{H^{1,p}(M)} := \left( \int_{M} |\nabla_g u|^p\; dv_g  + \int_{M} |u|^p\; dv_g \right)^{\frac{1}{p}}\; .
\]

The following definitions and notations related to (\ref{AB}) are
quite natural when one desires to introduce sharp $L^p$-entropy
inequalities in the Riemannian context.

The {\bf first best $L^p$-entropy constant} is defined by

\[
{\cal A}_0(p,g) := \inf \{ {\cal A} \in \R: \mbox{ there exists} \hspace{0,18cm} {\cal B} \in \R \hspace{0,18cm} \mbox{such that (\ref{AB})} \hspace{0,18cm} \mbox{holds for all} \hspace{0,18cm} u \in H^{1,p}(M) \hspace{0,18cm} \mbox{with} \hspace{0,18cm} \int_M |u|^p dv_g = 1\}\; .
\]

\n It follows directly that ${\cal A}_0(p,g)$ is well defined and moreover, by Jensen's inequality applied to (\ref{AB}), ${\cal A}_0(p,g)$ is positive for any $1 \leq p < n$.

The {\bf first sharp $L^p$-entropy inequality} states that there exists a constant ${\cal B} \in \R$ such that, for any $u \in H^{1,p}(M)$ with $\int_M |u|^p dv_g = 1$,

\[
\int_{M} |u|^p \log |u|^p\; dv_g \leq \frac{n}{p} \log \left( {\cal A}_0(p,g) \int_M |\nabla_g u|^p\; dv_g + {\cal B} \right)\; .
\]

If the preceding inequality is true, then we can define the {\bf second best $L^p$-entropy constant} as

\[
{\cal B}_0(p,g) := \inf \{ {\cal B} \in \R:\; (L({\cal A}_0(p,g),{\cal B})) \hspace{0,18cm} \mbox{holds for all} \hspace{0,18cm} u \in H^{1,p}(M) \hspace{0,18cm} \mbox{with} \hspace{0,18cm} \int_M |u|^p dv_g = 1\}
\]
and the {\bf second sharp $L^p$-entropy inequality} as the saturated version of ($L({\cal A}_0(p,g),{\cal B})$) on the functions $u \in H^{1,p}(M)$ with $\int_M |u|^p dv_g = 1$, that is

\[
\int_{M} |u|^p \log |u|^p\; dv_g \leq \frac{n}{p} \log \left( {\cal A}_0(p,g) \int_M |\nabla_g u|^p\; dv_g + {\cal B}_0(p,g) \right)\; .
\]

\n Note that ($L({\cal A}_0(p,g),{\cal B}_0(p,g))$) is sharp with respect to both the first and second best constants in the sense that none of them can be lowered.

In a natural way, one introduces the notion of extremal functions of ($L({\cal A}_0(p,g),{\cal B}_0(p,g))$). A function $u_0  \in H^{1,p}(M)$ satisfying $\int_M |u_0|^p dv_g = 1$ is said to be extremal, if

\[
\int_{M} |u_0|^p \log |u_0|^p\; dv_g = \frac{n}{p} \log \left( {\cal A}_0(p,g) \int_M |\nabla_g u_0|^p\; dv_g + {\cal B}_0(p,g) \right)\; .
\]

\n We denote by ${\cal E}_0(p,g)$ the set of all extremal functions of ($L({\cal A}_0(p,g),{\cal B}_0(p,g))$).

In \cite{Bro}, Brouttelande proved for dimensions $n \geq 3$ that ${\cal A}_0(2,g) = {\cal A}_0(2)$ and ($L({\cal A}_0(2,g),{\cal B})$) is valid for some constant ${\cal B} \in \R$. Subsequently, in \cite{Bro1}, he obtained the lower bound

\[
{\cal B}_0(2,g) \geq \frac{1}{2 n \pi e} \max_M R_g\; ,
\]

\n where $R_g$ stands for the scalar curvature of the metric $g$, and proved that if the inequality is strict, then the set ${\cal E}_0(2,g)$ is non-empty.

A simple lower bound for ${\cal B}_0(p,g)$ involving the volume $v_g(M)$ of $M$ also follows by taking a normalized constant function in ($L({\cal A}_0(p,g),{\cal B}_0(p,g))$), namely

\[
{\cal B}_0(p,g) \geq v_g(M)^{-p/n}\; ,
\]

\n provided that ${\cal B}_0(p,g)$ is well defined.

Our main contributions are gathered in the following result:

\begin{teor} \label{MT}
Let $(M,g)$ be a smooth closed Riemannian manifold of dimension $n \geq 2$. For any $1 < p \leq 2$ and $p < n$, we have:

\begin{itemize}

\item[{\bf (a)}] ${\cal A}_0(p,g) = {\cal A}_0(p)$;

\item[{\bf (b)}] there exists a constant ${\cal B} \in \R$ such that ($L({\cal A}_0(p,g),{\cal B})$) holds for all function $u \in H^{1,p}(M)$ with $\int_M |u|^p dv_g = 1$.

\end{itemize}

\end{teor}

Note that the above result extends to $1 < p \leq 2$ the corresponding one due to Brouttelande. However, his arguments do not apply to $p \neq 2$, so we here develop an alternative approach in order to prove Theorem \ref{MT}.

The idea of using subcritical interpolation inequalities for obtaining entropy inequalities was introduced by Bakry, Coulhon,
Ledoux and Sallof-Coste in \cite{Ba} within the Euclidean
environment. Namely, for $1 \leq p < n$, they showed how to produce
non-sharp $L^p$-entropy inequalities as a limit case of a class of
non-sharp Gagliardo-Nirenberg inequalities.

Later, this view point was explored in the sharp sense by Del Pino and Dolbeault \cite{DPDo} in order to establish the inequality
(\ref{dee}) for $1 < p < n$. Indeed, they considered a family of sharp Gagliardo-Nirenberg inequalities, interpolating
the $L^p$-Sobolev and $L^p$-entropy inequalities, whose extremal functions are explicitly known.

In trying to adapt the same idea of getting entropy inequalities as a limit case of subcritical interpolation inequalities to the
Riemannian context, the situation changes drastically mainly because extremal functions and second best constants are usually
unknown.

With the aim to make clear to the reader our program of proof and its main points of difficulty, a brief overview on related sharp Riemannian
Nash inequalities should be presented.

Let $1 < p \leq 2$ and $1 \leq q < p$. In \cite{CM2}, Ceccon and Montenegro established the existence of a constant $B \in \R$ such that the sharp $L^p$-Nash inequality

\begin{gather}\label{AB1}
\left( \int_M |u|^p\; dv_g \right)^{\frac{1}{\theta}} \leq \left( A_{opt} \int_M |\nabla_g u|^p\; dv_g + B \int_M |u|^p\; dv_g \right) \left(
\int_M |u|^q\; dv_g \right)^{\frac{p(1 - \theta)}{q \theta}} \tag{$I_{p,q}(A,B)$}
\end{gather}

\n holds for all function $u \in H^{1,p}(M)$, where $\theta = \frac{n(p - q)}{q(p - n) + np}$ and $A_{opt}$ is the first best $L^p$-Nash constant.

One also knows that

\begin{equation} \label{eq}
A_{opt} = N(p,q) \; ,
\end{equation}

\n where $N(p,q)$ stands for the best Euclidean $L^p$-Nash constant, see \cite{DHV} for $p = 2$ and \cite{CM2} for $1 < p < 2$. Namely, $N(p,q)$ is the best possible constant in the $L^p$-Nash inequality

\[
\left( \int_{\R^n} |u|^p\; dx \right)^{\frac{1}{\theta}} \leq A \left( \int_{\R^n} |\nabla u|^p\; dx \right) \left(
\int_{\R^n} |u|^q\; dx \right)^{\frac{p(1 - \theta)}{q \theta}}
\]

\n which holds for all function $C^\infty_0(\R^n)$, provided that $n \geq 2$ and $1 \leq q < p$. This last inequality was first proved by Nash in \cite{Na}. An alternative proof was given by Beckner in \cite{Be1}. For $p \neq 2$, we refer to the recent work \cite{Ce} by Ceccon.

One then defines the {\bf second best $L^p$-Nash constant} as

\[
B(p,q,g) := \inf \{ B \in \R:\; (I_{p,q}(N(p,q),B)) \hspace{0,18cm} \mbox{holds for all} \hspace{0,18cm} u \in H^{1,p}(M)\}\; .
\]

\n In this case, the {\bf second sharp $L^p$-Nash inequality} automatically holds on $H^{1,p}(M)$, namely

\[
\left( \int_M |u|^p\; dv_g \right)^{\frac{1}{\theta}} \leq \left( N(p,q) \int_M |\nabla_g u|^p\; dv_g + B(p,q,g) \int_M |u|^p\; dv_g \right) \left( \int_M |u|^q\; dv_g \right)^{\frac{p(1 - \theta)}{q \theta}}\; .
\]

Our strategy consists in rearranging the above inequality and applying logarithmic of both sides, so that

\begin{equation} \label{AE}
\frac{p}{\theta} \log{\frac{||u||_{L^p(M)}}{||u||_{L^q(M)}}} \leq \log{\frac{N(p,q) \int_M |\nabla_g u|^p\; dv_g + B(p,q,g) \int_M |u|^p\; dv_g} {\left( \int_M |u|^q\; dv_g \right)^{\frac{p}{q}}} }\, .
\end{equation}

\n and then letting the limit as $q \rightarrow p^-$.

The success of this plan will be result of the following three statements for $1 < p \leq 2$ and $p < n$:

\begin{itemize}

\item[{\bf (A)}] $N(p,q)$ converges to ${\cal A}_0(p)$ as $q \rightarrow p^-$

\item[{\bf (B)}] $B(p,q,g)$ is bounded for any $q < p$ close to $p$

\item[{\bf (C)}] ${\cal A}_0(p,g) \geq {\cal A}_0(p)$

\end{itemize}

\n In fact, assume for a moment that {\bf (A)}, {\bf (B)} and {\bf (C)} are true. Letting $q \rightarrow p^-$ in (\ref{AE}) and using {\bf (A)} and {\bf (B)}, after straightforward computations, one obtains the inequality ($L({\cal A}_0(p),{\cal B})$) for some constant ${\cal B}$. Indeed, as can easily be checked,
 
\[
\lim_{q \rightarrow p^-} \frac{p}{\theta} \log{\frac{||u||_{L^p(M)}}{||u||_{L^q(M)}}} =
\frac{p^3}{n} \lim_{q \rightarrow p^-} \frac{1}{p - q} \log{\frac{||u||_{L^p(M)}}{||u||_{L^q(M)}}} = \frac{p}{n} \int_M
|u|^p \log (|u|^p)\; dv_g
\]

\n for any $u \in H^{1,p}(M)$ satisfying $||u||_{L^p(M)} = 1$, see the proof of Proposition \ref{P.1} for a similar computation. Thus, one has ${\cal A}_0(p,g) \leq {\cal A}_0(p)$ and, thanks to {\bf (C)}, the conclusion of Theorem \ref{MT} follows.

We then describe some ideas involved in the proof of each one of the claims {\bf (A)}-{\bf (C)}. We begin by addressing {\bf (A)} and {\bf (C)} since their proofs are shortest.

\n {\bf On the proof of {\bf (A)}.} This claim is proved in Section 2 and uses Jensen's inequality and the fact to be proved that the best Nash constant $N(p,q)$ is increasing on $q$.

\n {\bf On the proof of {\bf (C)}.} This claim is proved in Section 3. Its proof is based on estimates of Gaussian bubbles. Precisely, we consider the following test function in ($L({\cal A}_0(p,g),{\cal B}_0(p,g))$):

\[
u_\varepsilon(exp_{x_0}(x)) = \eta(x) \varepsilon^{-\frac{n}{p}} u_0(\frac{x}{\varepsilon})
\]

\n defined locally around a point $x_0$ on $M$, where $\eta$ denotes a cutoff function supported in a ball centered at $0$. After using Cartan's expansion of the metric $g$ around $x_0$ and estimating each involved integral for $\varepsilon > 0$ small enough, the desired conclusion follows.

\n {\bf On the proof of {\bf (B)}.} This claim is proved in Section 4. Since its proof is rather long and technique, for a better understanding two important steps are presented under the form of lemma.

It suffices to prove the assertion for each sequence $(q_k) \subset (1, q)$ such that $q_k \rightarrow p$ as $k \rightarrow + \infty$. For such a sequence and each $k$, we consider the functional

\[
J_k(u) = \left( \int_M |\nabla_g u|^p\; dv_g + C_k \int_M |u|^p\; dv_g \right) \left( \int_M |u|^{q_k}\; dv_g \right)^{\frac{p(1 -
\theta_{q_k})}{q_k \theta_{q_k}}}
\]

\n on the set ${\cal H} = \{ u \in H^{1,p}(M):\; ||u||_{L^p(M)} = 1 \}$, where $C_k$ is defined as

\[
C_k = \frac{B(p,q_k,g) - (p-q_k)}{N(p,q_k)}\; .
\]

\n From its definition, we readily have $C_k < B(p,q_k,g) N(p,q_k)^{-1}$. Therefore, from the definition of $B(p,q_k,g)$, there exists a function $w_k \in {\cal H}$ satisfying $J_k(w_k) < N(p,q_k)^{-1}$, so that

\[
\inf_{u \in {\cal H}} J_k (u) < N(p,q_k)^{-1}\; .
\]

\n As usual, this last inequality leads to the existence of a $C^1$ minimizer $u_k \in {\cal H}$ for the functional $J_k$ on ${\cal H}$. The next steps consists in studying some fine properties satisfied by the sequence $(u_k)$ such as concentration and pointwise estimates. The main tools used here are blow-up method and Moser's iteration on elliptic PDEs.

We conclude the introduction by exposing some still open problems on sharp Riemannian entropy inequalities and related best constants.

Perhaps, contrary to what one might expect, it is not clear that the first best entropy constant ${\cal A}_0(p,g)$ is well defined and is equal to ${\cal A}_0(p)$ for all $p \geq n$. The great difficult in Riemannian $L^p$-entropy inequalities is that local-to-global type arguments do not usually work well. For example, the normalization condition $\int_M |u|^p dv_g = 1$ and the involved log functions in (\ref{AB}) do not allow a direct comparison to the corresponding flat case. In particular, it does not seem immediate that (\ref{AB}) is valid for some constants ${\cal A}$ and ${\cal B}$ and also that, for each $\varepsilon > 0$, there exists a constant ${\cal B}_\varepsilon \in \R$ such that

\[
\int_{M} |u|^p \log |u|^p\; dv_g \leq \frac{n}{p} \log \left( ({\cal A}_0(p) + \varepsilon) \int_M |\nabla_g u|^p dv_g + {\cal B}_\varepsilon \right)
\]

\n for all function $u \in H^{1,p}(M)$ with $\int_M |u|^p dv_g = 1$.

According to our contributions, the first best entropy constant ${\cal A}_0(p,g)$ is well defined for all $1 \leq p < n$ and equal to ${\cal A}_0(p)$ for all $1 < p \leq 2$ and $p < n$. In addition, based on developments over thirty years in the field of sharp entropy and Sobolev inequalities (see \cite{Au4}, \cite{DH} and references therein), one expects positive answers for the following questions:\\

\n {\bf Open problem 1.} Is ${\cal A}_0(p,g)$ well defined for all $p \geq n$?\\

\n {\bf Open problem 2.} Does ${\cal A}_0(p,g) = {\cal A}_0(p)$ hold in the three cases $p = 1,  n \geq 2$, $p = 2, n = 2$ or $p > 2, n \geq 2$?\\

\n {\bf Open problem 3.} Is ($L({\cal A}_0(p,g),{\cal B})$) valid for some constant ${\cal B}$ in the two cases $p = 1,  n \geq 2$ or $p = 2, n = 2$?\\

\n {\bf Open problem 4.} Is ($L({\cal A}_0(p,g),{\cal B})$) non-valid whenever $p > 2$,  $n \geq 2$ and $(M,g)$ has positive scalar curvature somewhere?

\section{Proof of the assertion {\bf (A)}}

This section is devoted to the proof of the proposition.

\begin{propo} \label{P.1}
Let $n \geq 2$ and $1 \leq q < p$. We have

\begin{equation} \label{LN}
\lim_{q \rightarrow p^-} N(p,q) = {\cal A}_0(p) \; .
\end{equation}
\end{propo}

We first prove that

\begin{equation} \label{DN}
N(p,q) \leq {\cal A}_0(p)
\end{equation}

\n for all $1 < q < p$.

\n Let $u \in C^{\infty}_0(\R^n)$ such that $||u||_{L^p(\R^n)} = 1$. Using Jensen's inequality, we have

\[
-\ln \left(\int_{\R^n} |u|^q\; dx \right) = -\ln \left(\int_{\R^n} |u|^{q - p} |u|^p\; dx \right) \leq - \int_{\R^n} \ln(|u|^{q - p})
|u|^p\; dx = \frac{p - q}{p} \int_{\R^n} \ln(|u|^p) |u|^p\; dx\; .
\]

\n Joining this inequality with (\ref{dee}), one obtains

\[
\left(\int_{\R^n} |u|^q dx \right)^{- \frac{p^2}{n(p-q)}} \leq
{\cal A}_0(p) \int_{\R^n} |\nabla u|^p dx \; .
\]

\n Using the fact that $||u||_{L^p(\R^n)} = 1$, we then derive the Nash inequality

\[
\left( \int_{\R^n} |u|^p dx \right)^{\frac{1}{\theta}} \leq {\cal A}_0(p)
\left( \int_{\R^n} |\nabla u|^p dx \right) \left(\int_{\R^n} |u|^q dx
\right)^{\frac{p(1 - \theta)}{\theta q}}\; ,
\]

\n where

\[
\theta = \frac{n(p - q)}{qp - qn + np} \; .
\]

\n By homogeneity, the above inequality is valid for all function $u \in C^{\infty}_0(\R^n)$, so that the assertion (\ref{DN}) follows.

We now prove that $N(p,q)$ is monotonically increasing on $q$.

\n Let $1 < q_1 < q_2 < p$ fixed. An usual interpolation inequality yields

\[
\left( \int_{\R^n} |u|^{q_2}\; dx \right)^{\frac{1}{q_2}} \leq \left(\int_{\R^n} |u|^{q_1}\; dx \right)^{\frac{\mu}{q_1}}
\left(\int_{\R^n} |u|^{p}\; dx \right)^{\frac{1 - \mu}{p}}\; ,
\]

\n where $\frac{1}{q_2} = \frac{\mu}{q_1} + \frac{1 - \mu}{p}$.

\n Plugging this inequality in

\[
\left(\int_{\R^n} |u|^{p}\; dx \right)^{\frac{1}{\theta_2}} \leq N(p,q_2) \left( \int_{\R^n} |\nabla u|^p\; dx \right)
\left( \int_{\R^n} |u|^{q_2}\; dx \right)^{\frac{p(1 - \theta_2)}{q_2 \theta_2}} \; ,
\]

\n where $\theta_2 = \frac{n(p - q_2)}{q_2(p - n) + np}$, one gets

\[
\left(\int_{\R^n} |u|^p\; dx \right)^{1 + \mu\frac{1 - \theta_2}{\theta_2}} \leq N(p,q_2)
\left(\int_{\R^n} |\nabla u|^p\; dx \right) \left( \int_{\R^n} |u|^{q_1}\; dx\right)^{\frac{p}{q_1} (\mu
\frac{1 - \theta_2}{\theta_2})}\; .
\]

\n On the other hand, the definition of $\mu$ produces the relations

\[
1 + \mu \frac{1 -
\theta_2}{\theta_2} = \frac{1}{\theta_1}
\]

\n and

\[
\frac{p}{q_1}\left( \mu \frac{1 -
\theta_2}{\theta_2}\right) =
\frac{p}{q_1}\frac{1 - \theta_1}{\theta_1} \; ,
\]

\n where $\theta_1 = \frac{n(p - q_1)}{q_1(p - n) + np}$, so that $N(p,q_1) \leq N(p,q_2)$ and the monotonicity follows.

\n We then denote

\begin{equation}
A(p) = \lim_{q \rightarrow p^-} N(p,q)\; .
\end{equation}

\n It is clear by (\ref{DN}) that $A(p) \leq {\cal A}_0(p)$.

The remaining of the proof is devoted to show that $A(p) \geq {\cal A}_0(p)$.

\n Rearranging the sharp Nash inequality and taking logarithm of both sides, one has

\[
\frac{1}{\theta} \log \left(\frac{||u||_p}{||u||_q} \right) \leq \log \left(N(p,q) \frac{||\nabla u||_p^p}{||u||_q^p}\right)^{\frac{1}{p}}\; .
\]

\n Using the definition of $\theta$ and taking the limit on $q$, one gets

\[
\frac{p^2}{n} \lim_{q \rightarrow p^-} \frac{1}{p - q} \log \left(\frac{||u||_p}{||u||_q} \right)\leq \log \left(A(p) \frac{||\nabla u||_p^p}{||u||_p^p}\right)^{\frac{1}{p}}\; .
\]

\n We now compute the left-hand side limit. We first write

\[
\log \left(\frac{||u||_p}{||u||_q} \right) = \frac{1}{p} \log(||u||_p^p) - \frac{1}{q} \log(||u||_q^q) = \frac{q - p}{p}
\log(||u||_q) + \frac{1}{p} \left( \log(||u||_p^p) - \log(||u||_q^q) \right)\; .
\]

\n A straightforward computation then gives

\[
\lim_{q \rightarrow p^-} \frac{1}{p - q} \log \left(\frac{||u||_p}{||u||_q} \right) = \frac{1}{p} \int_{\R^n}
\frac{|u|^p}{||u||_p^p} \log \left(\frac{|u|}{||u||_p}\right)\; dx\; .
\]

\n Therefore,

\[
\int_{\R^n} \frac{|u|^p}{||u||_p^p} \log \left( \frac{|u|^p}{||u||_p^p} \right)\; dx \leq \frac{n}{p} \log \left( A(p) \frac{||\nabla
u||_p^p}{||u||_p^p}\right)\; ,
\]

\n or equivalently,

\[
\int_{\R^n}|u|^p \log (|u|^p)\; dx \leq \frac{n}{p} \log \left( A(p)
\int_{\R^n} |\nabla u|^p\; dx \right)
\]

\n for all function $u \in C^\infty_0(\R^n)$ with $||u||_{L^p(\R^n)} = 1$. So, $A(p) \geq {\cal A}_0(p)$ and the proof of Proposition \ref{P.1} follows. \bl

\section{Proof of the assertion {\bf (C)}}

In this section, we prove the following result:

\begin{propo}\label{P.2}
For each $n \geq 2$ and $1 < p < n$, we have

\[
{\cal A}_0(p,g) \geq {\cal A}_0(p)\; .
\]

\end{propo}

Using the assumption that $n \geq 2$ and $1 < p < n$, one gets constants ${\cal A}$ and ${\cal B}$ such that $L({\cal A},{\cal B})$ is valid. It suffices to show that ${\cal A} \geq {\cal A}_0(p)$.

One knows that $L({\cal A},{\cal B})$ is equivalent to

\begin{equation} \label{EF1}
\frac{1}{||u||_p^p} \int_M |u|^p \log(|u|^p)\; dv_g + (\frac{n}{p} - 1) \log (||u||_p^p) \leq \frac{n}{p} \log \left(  {\cal A} \int_M |\nabla
u|_g^p\; dv_g + {\cal B} \int_M |u|^p\; dv_g \right)
\end{equation}

\n for all function $u \in C^{\infty}(M)$.

We first fix a point $x_0 \in M$ and an extremal function $u_0(x) = a e^{-b|x|^{\frac{p}{p-1}}} \in W^{1,p}(\R^n)$ for the sharp entropy inequality (\ref{dee}) (see \cite{DPDo} or \cite{G}), where $a$ and $b$ are positive constants chosen so that $||u_0||_{L^p(\R^n)} = 1$.

Consider now a geodesic ball $B(x_0,\delta) \subset M$ and a radial cutoff function $\eta \in C^\infty(B(0,\delta))$ satisfying $\eta = 1$ in $B(0,\frac{\delta}{2})$, $\eta = 0$ outside $B(0,\delta),\ 0 \leq \eta \leq 1$ in $B(0,\delta)$. For $\varepsilon > 0$ and $x \in B(0,\delta)$, set

\[
u_\varepsilon(exp_{x_0}(x)) = \eta(x) \varepsilon^{-\frac{n}{p}} u_0(\frac{x}{\varepsilon})\; .
\]

The asymptotic behavior of each integral of (\ref{EF1}) computed at $u_\varepsilon$ with $\varepsilon > 0$ small enough is now presented.

\n Denote

\[
I_1 = \int_{\R^n} u_0^p \log(u_0^p)\; dx, \ \ I_2 = \int_{\R^n} |\nabla u_0|^p\; dx\; ,
\]

\[
J_1 = \int_{\R^n} u_0^p |x|^2\; dx, \ \ J_2 = \int_{\R^n} |\nabla u_0|^p |x|^2\; dx,\ \ J_3 = \int_{\R^n} u_0^p \log(u_0^p) |x|^2\; dx\; .
\]

\n Using the expansion of volume element in geodesic coordinates

\[
\sqrt{det g} = 1 - \frac{1}{6} \sum_{i,j = 1}^{n} Ric_{ij}(x_0) x_i x_j + O(r^3)\; ,
\]

\n where $Ric_{ij}$ denotes the components of the Ricci tensor in these coordinates and $r=|x|$, one easily checks that

\begin{equation} \label{1}
\int_{M} u_\varepsilon^p\; dv_g = 1 - \frac{R_g(x_0)}{6n} J_1 \varepsilon^2 + O(\varepsilon^4)\; ,
\end{equation}

\begin{equation} \label{2}
\int_{M} u_\varepsilon^p \log(u_\varepsilon^p)\; dv_g = I_1 - n \log \varepsilon - \frac{R_g(x_0)}{6n} J_3 \varepsilon^2 + \frac{R_g(x_0)}{6} J_1 \varepsilon^2 \log \varepsilon  + o(\varepsilon^4 \log \varepsilon)
\end{equation}

\n and

\begin{equation} \label{3}
\int_{M} |\nabla u_\varepsilon|^p\; dv_g = \varepsilon^{-p} \left(I_2 - \frac{R_g(x_0)}{6n} J_2 \varepsilon^2 + O(\varepsilon^4) \right)\; .
\end{equation}

\n Plugging $u_\varepsilon$ in (\ref{EF1}), from (\ref{1}), (\ref{2}) and (\ref{3}), one obtains

\[
\frac{I_1 - n \log \varepsilon - \frac{R_g(x_0)}{6n} J_3 \varepsilon^2 + \frac{R_g(x_0)}{6} J_1 \varepsilon^2 \log \varepsilon  + o(\varepsilon^4 \log \varepsilon)}{1 - \frac{R_g(x_0)}{6n} J_1 \varepsilon^2 + O(\varepsilon^4)} + (\frac{n}{p} - 1) \log \left(1 - \frac{R_g(x_0)}{6n} J_1 \varepsilon^2 + O(\varepsilon^4)\right)
\]

\[
\leq -n \log \varepsilon + \frac{n}{p} \log\left( {\cal A} I_2 - \frac{R_g(x_0)}{6n} {\cal A} J_2 \varepsilon^2 + {\cal B} \varepsilon^p +
O(\varepsilon^q)\right)
\]

\n for $\varepsilon > 0$ small enough, where $q = \min\{4, p +2\}$.

\n Taylor`s expansion then guarantees that

\begin{equation}\label{taylor}
I_1 - n \log \varepsilon - \frac{R_g(x_0)}{6n} J_3 \varepsilon^2 +
\frac{R_g(x_0)}{6} J_1 \varepsilon^2 \log \varepsilon +
\frac{R_g(x_0)}{6n} I_1 J_1 \varepsilon^2 - \frac{R_g(x_0)}{6} J_1
\varepsilon^2 \log \varepsilon
\end{equation}

\[
- (\frac{n}{p} - 1) \frac{R_g(x_0)}{6n} J_1 \varepsilon^2 +
O(\varepsilon^4 \log \varepsilon) \leq - n \log \varepsilon +
\frac{n}{p} \log ({\cal A} \; I_2) - \frac{n}{p} \frac{R_g(x_0)}{6n}
\frac{J_2}{I_2} \varepsilon^2 + \frac{n}{p} \frac{{\cal B}}{{\cal A} \; I_2}
\varepsilon^p + O(\varepsilon^q)\; .
\]

\n After a suitable simplification, one arrives at

\[
I_1 - \frac{n}{p} \log({\cal A} \; I_2) \leq - \frac{R_g(x_0)}{6n} \left( I_1 J_1 + \frac{n}{p} \frac{J_2}{I_2} - (\frac{n}{p} - 1) J_1 -
J_3\right) \varepsilon^2 + O(\varepsilon^p)
\]

\n for $\varepsilon > 0$ small enough, so that

\[
I_1 \leq \frac{n}{p} \log({\cal A} \; I_2) \; .
\]

\n Thus, since $u_0$ is an extremal function of (\ref{dee}), one has

\[
I_1 = \frac{n}{p} \log({\cal A}_0(p) \; I_2)\; ,
\]

\n so that ${\cal A} \geq {\cal A}_0(p)$, and the conclusion of Proposition \ref{P.2} readily follows. \bl

\section{Proof of the assertion {\bf (B)}}

This section is devoted to the proof of the following theorem:

\begin{teor}\label{lB}
Let $(M,g)$ be a smooth closed Riemannian manifold of dimension $n \geq 2$. For each fixed $1 < p \leq 2$ and $p < n$, the best constant $B(p,q,g)$ is bounded for any $q < p$ close to $p$.
\end{teor}

Clearly, it suffices to prove this result for an arbitrary sequence $(q_k) \subset (1, p)$ converging to $p$ as $k \rightarrow + \infty$.

\n Define

\[
C_k = \frac{B(p,q_k,g) - (p-q_k)}{N(p,q_k)}\; .
\]

\n Since $C_k < B(p,q_k,g) N(p,q_k)^{-1}$, according to the definition of the best constant $B(p,q_k,g)$, we have

\begin{equation}\label{ldf}
\inf_{u \in {\cal H}} J_k (u) < N(p,q_k)^{-1}\; ,
\end{equation}

\n where ${\cal H} = \{ u \in H^{1,p}(M):\; ||u||_{L^p(M)} = 1 \}$ and

\[
J_k(u) = \left( \int_M |\nabla_g u|^p\; dv_g + C_k \int_M |u|^p\; dv_g \right) \left( \int_M |u|^{q_k}\; dv_g \right)^{\frac{p(1 -
\theta_k)}{q_k \theta_k}}\; .
\]

\n As can easily be checked, the functional $J_k$ is of class $C^1$.

\n The condition (\ref{ldf}) and an usual argument of direct minimization lead to a minimizer $u_k \in {\cal H}$ of $J_k$, so that

\begin{equation}\label{l3iha}
\nu_k := J_k(u_k) = \inf_{u \in {\cal H}} J_k(u)\; .
\end{equation}

\n Note that we can assume $u_k \geq 0$, since $|\nabla_g |u_k|| = |\nabla_g u_k|$ almost everywhere.

\n In addition, since $J_k$ is differentiable, $u_k$ satisfies the quasilinear elliptic equation

\begin{equation}\label{l3ep}
A_k \Delta_{p,g} u_k + A_k C_k u_k^{p - 1} + \frac{1 -
\theta_k}{\theta_k} B_k u_k^{q_k - 1} = \frac{\nu_k}{\theta_k} u_k^{p - 1}
\hspace{0,2cm} \mbox{on} \hspace{0,2cm} M\; ,
\end{equation}

\n where $\Delta_{p,g} = -{\rm div}_g(|\nabla_g|^{p-2} \nabla_g)$ is the $p$-Laplace operator of $g$,

\[
\theta_k = \frac{n(p - q_k)}{np + pq_k - nq_k}\; ,
\]

\[
A_k = \left(\int_M u_k^{q_k}\; dv_g \right)^{\frac{p( 1 - \theta_k)}{q_k \theta_k}}\; ,
\]

\[
B_k = \left(\int_M |\nabla_g u_k|^p\; dv_g + C_k \right) \left(\int_M u_k ^{q_k}\; dv_g \right)^{\frac{p(1 - \theta_k)}{q_k \theta_k} - 1}\; .
\]

\n So, by Serrin \cite{Se}, $u_k \in L^\infty(M)$ and, by Tolksdorf's \cite{To}, it follows that $u_k$ is of class $C^1$.

\n A relation that will be useful is

\begin{equation} \label{nu}
B_k \int_M u_k^{q_k}\; dv_g =\nu_k\; .
\end{equation}

\n Modulo a subsequence, we analyze two possible situations for the sequence $(A_k)$.

\begin{itemize}
\item[{\bf (i)}] $\lim_{k \rightarrow + \infty} A_k > 0$;

\item[{\bf (ii)}] $\lim_{k \rightarrow + \infty} A_k = 0$.
\end{itemize}

\n Assume that {\bf (i)} occurs. Taking $u_k$ as a test function in (\ref{l3ep}) and after using Proposition \ref{P.1}, one gets

\begin{equation} \label{lim1}
A_k C_k \leq \nu_k \leq c
\end{equation}

\n for $k$ large enough and some constant $c$ independent of $k$, so that $(C_k)$ is bounded. Thus, the conclusion of Theorem \ref{lB} follows directly from the definition of $C_k$.

The remaining of this section is dedicated to prove the boundedness of $(C_k)$ by assuming that the situation {\bf (ii)} occurs. In this case, the inequality

\begin{equation} \label{lim2}
||u_k||_{L^\infty(M)} A_k^{\frac{n}{p^2}} \geq 1
\end{equation}

\n follows from

\[
1 = ||u_k||_{L^p(M)}^p = \int_M u_k^p\; dv_g \leq ||u_k||_{L^\infty(M)}^{p - q_k} \int_M u_k^{q_k}\; dv_g
\]

\n and implies that $(u_k)$ blows up in $L^\infty(M)$. Part of the proof consists in performing a fine analysis of concentration of the sequence $(u_k)$.

\n At first, we have

\begin{equation} \label{asym}
\lim \limits_{k \rightarrow + \infty} \nu_k = {\cal A}_0(p)^{-1}\; .
\end{equation}

\n In fact, a combination between the Nash inequality

\[
\left( \int_M |u|^p\; dv_g \right)^{\frac{1}{\theta_k}} \leq \left( N(p,q_k) \int_M |\nabla_g u|^p\; dv_g + B(p,q_k,g) \int_M |u|^p\; dv_g \right) \left(
\int_M |u|^{q_k}\; dv_g \right)^{\frac{p(1 - \theta_k)}{q_k \theta_k}}
\]

\n and the definition of $C_k$ yields

\[
N(p,q_k)^{-1} \leq \left(\int_M |\nabla_g u_k|^p\; dv_g + C_k \right)\left(\int_M u_k^{q_k}\; dv_g \right)^{\frac{p(1 - \theta_k)}{q_k \theta_k}} + \frac{(p - q_k)}{N(p,q_k)} \; A_k  = \nu_k + A_k \;\frac{p - q_k}{N(p,q_k)}\; .
\]

\n By Proposition \ref{P.1}, $N(p,q_k)$ remains away from zero for $k$ large enough. Therefore,

\[
\liminf \limits_{k \rightarrow + \infty} \nu_k \geq {\cal A}_0(p)^{-1}
\]

\n and so the limit (\ref{asym}) follows from (\ref{ldf}).

Let $x_k \in M$ be a maximum point of $u_k$, that is

\begin{equation}\label{l3ix}
u_k(x_k) = ||u_k||_{L^\infty(M)}\; .
\end{equation}

The following concentration property satisfied by the sequence $(u_k)$ plays an essential role in what follows. The main tools used in its proof are the blow-up method and Moser's iteration.

\begin{lema}\label{CA}
We have

\[
\lim \limits_{\sigma\rightarrow +\infty} \lim \limits_{k \rightarrow + \infty} \int_{B(x_k,\sigma A_k^{\frac{1}{p}})} u_k^p\; dv_g = 1\; .
\]
\end{lema}

\n {\bf Proof of Lemma \ref{CA}.} Let $\sigma > 0$. For each $x \in B(0, \sigma)$ and $k$ large, we define

\[
\begin{array}{l}
h_k(x) = g(\exp_{x_k} (A_k^{\frac{1}{p}} x))\; , \vspace{0,3cm}\\
\varphi_k(x) = A_k ^{\frac{n}{p^2}}
u_k(\exp_{x_k}(A_k^{\frac{1}{p}} x))\; .
\end{array}
\]

\n Note that the above expressions are well defined because we are assuming that $A_k \rightarrow 0$ as $k \rightarrow + \infty$ (situation {\bf (ii)}).

\n By (\ref{l3ep}) and (\ref{nu}), one easily deduces that

\[
\Delta_{p,h_k} \varphi_k + A_k C_k \varphi_k^{p - 1} + \frac{1 - \theta_k}{\theta_k} \nu_k \varphi_k^{q_k - 1} = \frac{\nu_k}{\theta_k} \varphi_k^{p - 1} \hspace{0,2cm} \mbox{on} \hspace{0,2cm}
B(0,\sigma) \; .
\]

\n Using the mean value theorem and the value of $\theta_k$, one gets

\begin{equation}\label{el1}
\Delta_{p,h_k} \varphi_k + A_k C_k \varphi_k^{p - 1} = \nu_k \left( \varphi_k^{q_k - 1} + \frac{pq_k - nq_k + np}{n} \varphi_k^{\rho_k} \log \varphi_k \right) \hspace{0,2cm} \mbox{on} \hspace{0,2cm} B(0,\sigma)
\end{equation}

\n for some $\rho_k \in (q_k - 1,p - 1)$.

\n Consider $\varepsilon > 0$ fixed such that $p + \varepsilon < \frac{np}{n - p}$. Since $\varphi_k^{\rho_k} \log (\varphi_k^\varepsilon) \leq \varphi_k^{p - 1 + \varepsilon}$, we then get

\[
\Delta_{p,h_k} \varphi_k + A_k C_k \varphi_k^{p - 1} \leq \nu_k \left( \varphi_k^{q_k - 1} + \frac{pq_k - nq_k + np}{n \varepsilon}
\varphi_k^{\varepsilon + p - 1} \right) \hspace{0,2cm} \mbox{on} \hspace{0,2cm} B(0,\sigma)\; .
\]

\n Once all coefficients of this equation are bounded, the Moser's iterative scheme (see \cite{Se}) produces

\[
(A_k^{\frac{n}{p^2}} ||u_k||_{L^{\infty}(M)})^p = \sup_{B(0,\sigma)} \varphi_k^p \leq c_\sigma \int_{B(0,2\sigma)} \varphi_k^p\; dv_{h_k} =
c_\sigma \int_{B(x_k,2 \sigma A_k^{\frac{1}{p}})} u_k^p\; dv_g \leq c_\sigma
\]

\n for all $k$ large enough, where $c_\sigma$ is a constant independent of $q$.

\n So, using (\ref{lim2}) in the above inequality, one readily obtains

\begin{equation}\label{lep}
1 \leq ||\varphi_k||_{L^{\infty}(B(0,\sigma))} \leq c_\sigma^{1/p}
\end{equation}

\n for all $k$ large enough.

\n By (\ref{lim1}), up to a subsequence, we have

\[
\lim_{k \rightarrow + \infty} A_k C_k = C \geq 0\; .
\]

\n Applying the Tolksdorf's elliptic theory to (\ref{el1}), thanks to (\ref{lep}), one easily checks that $\varphi_k \rightarrow \varphi$ in $C^1_{loc}(\R^n)$ and $\varphi \not \equiv 0$.

\n Letting now $k \rightarrow + \infty$ in (\ref{el1}) and using (\ref{asym}), one gets

\begin{equation} \label{ELim}
\Delta_{p,\xi} \varphi + C \varphi^{p - 1} = {\cal A}_0(p)^{-1} \left( \varphi^{p - 1} + \frac{p}{n} \varphi^{p - 1} \log (\varphi^p)\right) \hspace{0,2cm} \mbox{on} \hspace{0,2cm} \R^n \;,
\end{equation}

\n where $\Delta_{p,\xi}$ stands for the Euclidean $p$-Laplace operator.

\n For each $\sigma > 0$, we have

\begin{equation} \label{whole}
\int_{B(0,\sigma)} \varphi^p\; dx = \lim_{k \rightarrow + \infty} \int_{B(0,\sigma)} \varphi_k^p\; dv_{h_k} = \lim_{k \rightarrow + \infty}\int_{B(x_k,\sigma A_k^{\frac{1}{p}})} u_k^p\; dv_g \leq 1
\end{equation}

\n and, by $(\ref{ldf})$ and Proposition \ref{P.1},

\begin{equation} \label{grad}
\int_{B(0,\sigma)} |\nabla \varphi|^p\; dx = \lim_{k \rightarrow + \infty} \int_{B(0,\sigma)} |\nabla _{h_k} \varphi_k|^p\; dv_{h_k} = \lim_{k \rightarrow + \infty} \left( A_k \int_{B(x_k,\sigma A_k^{\frac{1}{p}})} |\nabla_g u_k|^p\; dv_g \right) \leq {\cal A}_0(p)^{-1}\; .
\end{equation}

\n In particular, one has $\varphi \in W^{1,p}(\R^n)$.

\n Consider then a sequence of nonnegative functions $(\varphi_k) \subset C^\infty_0(\R^n)$ converging to $\varphi$ in $W^{1,p}(\R^n)$. Taking $\varphi_k$ as a test function in (\ref{ELim}), we can write

\[
\frac{n}{p} {\cal A}_0(p) \int_{\R^n} |\nabla \varphi|^{p-2} \nabla \varphi \cdot \nabla \varphi_k\; dx + \frac{n}{p} {\cal A}_0(p) C \int_{\R^n} \varphi^{p-1} \varphi_k\; dx = \int_{\R^n} \varphi^{p-1} \varphi_k \log (\varphi^p)\; dx + \frac{n}{p} \int_{\R^n} \varphi^{p-1} \varphi_k \; dx
\]

\[
= \int_{[\varphi \leq 1]} \varphi^{p-1} \varphi_k \log (\varphi^p)\; dx + \int_{[\varphi \geq 1]} \varphi^{p-1} \varphi_k \log (\varphi^p)\; dx + \frac{n}{p} \int_{\R^n} \varphi^{p-1} \varphi_k \; dx
\]

\n Letting $k \rightarrow + \infty$ and applying Fatou's lemma and the dominated convergence theorem in the right-hand side, one gets

\[
\frac{n}{p} {\cal A}_0(p) \int_{\R^n} |\nabla \varphi|^p\; dx \leq \int_{[\varphi \leq 1]} \varphi^p \log (\varphi^p)\; dx + \int_{[\varphi \geq 1]} \varphi^p \log (\varphi^p)\; dx + \frac{n}{p} \int_{\R^n} \varphi^p\; dx
\]

\[
= \int_{\R^n} \varphi^p \log (\varphi^p)\; dx + \frac{n}{p} \int_{\R^n} \varphi^p\; dx\; .
\]

\n Rewriting this inequality in function of $\psi(x) = \frac{\varphi(x)}{||\varphi||_p}$, one has

\[
\frac{n}{p} {\cal A}_0(p) \int_{\R^n} |\nabla \psi|^p\; dx \leq \int_{\R^n} \psi^p \log (\psi^p)\; dx + \frac{n}{p} \int_{\R^n} \psi^p\; dx + \log ||\varphi||_p^p\; .
\]

\n Note that this inequality combined with

\[
\int_{\R^n} \psi^p \log (\psi^p)\; dx  \leq \frac{n}{p} \log \left({\cal A}_0(p) \int_{\R^n} |\nabla \psi|^p\; dx \right)
\]

\n produces

\[
{\cal A}_0(p) \int_{\R^n} |\nabla \psi|^p\; dx \leq \log \left({\cal A}_0(p)\; \int_{\R^n} |\nabla \psi|^p\; dx\right) + 1 + \frac{p}{n} \log ||\varphi||_p^p\; .
\]

\n Since $\log x \leq x - 1$ for all $x > 0$, we find $\log ||\varphi||_{L^p(\R^n)}^p \geq 0$ or, equivalently, $||\varphi||_{L^p(\R^n)} \geq 1$.

\n Thus, by (\ref{whole}), we conclude that $||\varphi||_{L^p(\R^n)} = 1$, so that

\[
\lim_{\sigma \rightarrow \infty} \lim_{k \rightarrow + \infty} \int_{B(x_k,\sigma A_k^{\frac{1}{p}})} u_k^p\; dv_g = \int_{\R^n} \varphi^p\; dx = 1\; .
\]
\bl

By using the concentration property provided in Lemma \ref{CA}, we now establish a pointwise estimate for the sequence $(u_k)$. This result is the key ingredient in the final part of the proof of Theorem \ref{lB}. The necessary tools in its proof are the same ones of the previous proof.

\begin{lema}\label{UE}
For any $\lambda > 0$, there exists a constant $c_\lambda > 0$, independent of $q < p$, such that

\[
d_g(x,x_k)^\lambda u_k(x) \leq c_\lambda A_k^{\frac{\lambda}{p} - \frac{n}{p^2}}
\]

\n for all $x \in M$ and $k$ large enough, where $d_g$ stands for the distance with respect to the metric $g$.
\end{lema}

\n {\bf Proof of Lemma \ref{UE}.} Suppose by contradiction that the above statement fails. Then, there exist $\lambda_0 > 0$ and $y_k \in M$ for each $k$ such that $f_{k,\lambda_0}(y_k) \rightarrow + \infty$ as $k \rightarrow + \infty$, where

\[
f_{k,\lambda}(x) = d_g(x,x_k)^\lambda u_k(x) A_k^{\frac{n}{p^2} - \frac{\lambda}{p}} \; .
\]

\n Without loss of generality, we assume that $f_{k,\lambda_0}(y_k) = ||f_{k,\lambda_0}||_{L^\infty(M)}$.

\n From (\ref{lep}), we have

\[
f_{k,\lambda_0}(y_k) \leq c \frac{u_k(y_k)}{||u_k||_{\infty}} d_g(x_k,y_k)^{\lambda_0} A_k^{- \frac{\lambda_0}{p}} \leq c d_g(x_k,y_k)^{\lambda_0} A_k^{- \frac{\lambda_0}{p}}\; ,
\]

\n so that

\begin{equation}\label{inf}
d_g(x_k,y_k) A_k^{- \frac{1}{p}} \rightarrow + \infty\ \ {\rm as}\ \ k \rightarrow + \infty\; .
\end{equation}

\n For any fixed $\varepsilon \in (0,1)$ and $\sigma > 0$, we first claim that

\begin{equation}\label{3int}
B(y_k,\varepsilon d(x_k,y_k)) \cap B(x_k, \sigma A_k^{\frac{1}{p}}) = \emptyset
\end{equation}

\n for $k$ large enough.

\n Clearly, this assertion follows from

\[
d_g(x_k,y_k) \geq \sigma  A_k^{\frac{1}{p}} + \varepsilon d(x_k,y_k)
\]

\n or equivalently,
\[
d_g(x_k, y_k)(1 - \varepsilon) A_k^{- \frac{1}{p}} \geq \sigma \; .
\]

\n But the above inequality is automatically satisfied, since $d_g(x_k,y_k) A_k^{- \frac{1}{p}} \rightarrow + \infty$ as $k \rightarrow + \infty$.

\n We assert that exists a constant $c > 0$, independent of $k$, such that

\begin{equation}\label{lim}
u_k(x) \leq c u_k(y_k)
\end{equation}

\n for all $x \in B(y_k, \varepsilon d_g(x_k,y_k))$ and $k$ large enough. Indeed,

\[
d_g(x,x_k) \geq d_g(x_k,y_k) - d_g(x,y_k) \geq (1 - \varepsilon) d_g(x_k,y_k)
\]

\n for all $x \in B(y_k, \varepsilon d_g(x_k,y_k))$. Thus,

\[
d_g(x_k,y_k)^{\lambda_0} u_k(y_k) A_k^{\frac{n}{p^2} - \frac{\lambda_0}{p}} = f_{k,\lambda_0}(y_k) \geq f_{k,\lambda_0}(x) = d_g(x,x_k)^{\lambda_0} u_k(x) A_k^{\frac{n}{p^2} -
\frac{\lambda_0}{p}}
\]

\[
\geq (1 - \varepsilon)^{\lambda_0} d_g(x_k,y_k)^{\lambda_0} u_k(x) \; A_k^{\frac{n}{p^2} - \frac{\lambda_0}{p}} \; ,
\]

\n so that

\[
u_k(x) \leq \left(\frac{1}{1 - \varepsilon}\right)^{\lambda_0} u_k(y_k)
\]

\n for all $x \in B(y_k, \varepsilon d_g(x_k,y_k))$ and $k$ large enough, as claimed.

\n We now define

\[
\begin{array}{l}
\tilde{h}_k(x) = g(\exp_{y_k} (A_k^{\frac{1}{p}} x))\; , \vspace{0,3cm}\\
\tilde{\varphi}_k(x) = A_k ^{\frac{n}{p^2}} u_k(\exp_{y_k}(A_k^{\frac{1}{p}} x)) \; .
\end{array}
\]

\n By (\ref{l3ep}) and (\ref{nu}), it readily follows that

\[
\Delta_{p,h_k} \tilde{\varphi}_k + A_k C_k \tilde{\varphi}_k^{p - 1} + \frac{1 - \theta_k}{\theta_k} \nu_k \; \tilde{\varphi}_k^{q_k - 1}
= \frac{\nu_k}{\theta_k} \tilde{\varphi}_k^{p - 1} \hspace{0,2cm} \mbox{on} \hspace{0,2cm} B(0,3) \; .
\]

\n Applying the mean value theorem, one obtains

\begin{equation}\label{leqgn}
\Delta_{p,h_k} \tilde{\varphi}_k + A_k C_k \tilde{\varphi}_k^{p - 1} = \nu_k \left( \tilde{\varphi}_k^{q_k - 1} + \frac{pq_k - nq_k +
np}{n} \tilde{\varphi}_k^{\rho_k} \log (\tilde{\varphi}_k)\right) \hspace{0,2cm} \mbox{on} \hspace{0,2cm} B(0,3) \; ,
\end{equation}

\n where $\rho_k \in (q_k - 1,p - 1)$.

\n For fixed $\varepsilon > 0$ such that $p + \varepsilon < \frac{np}{n - p}$, one has $\tilde{\varphi}_k^{\rho_k} \log (\tilde{\varphi}_k^\varepsilon)
\leq \tilde{\varphi}_k^{p - 1 + \varepsilon}$. So, the Moser's iterative scheme applied to (\ref{leqgn}) yields

\begin{equation}\label{pf}
\mu_k^{\frac{p}{p - q_k}} = (A_k^{\frac{n}{p^2}} u_k(y_k))^p \leq \sup_{B(0,1)} \tilde{\varphi}_k^p \leq c \int_{B(0,2)}
\tilde{\varphi}_k^p\; dv_{\tilde{h}_k} = c \int_{B(y_k,A_k^{\frac{1}{p}})} u_k^p\; dv_g
\end{equation}

\n for $k$ large enough, where

\begin{equation}\label{mu}
\mu_k = u_k(y_k)^{p - q_k} \int_M u_k^{q_k} dv_g\; .
\end{equation}

\n We next analyze two independent situations that can occur:

\begin{itemize}
\item[{\bf (I)}] $\mu_k \geq 1 - \theta_k$ for all $k$, up to a subsequence;

\item[{\bf (II)}] $\mu_k < 1 - \theta_k$ for $k$ large enough.
\end{itemize}

\n In each case, we derive a contradiction. If the assertion {\bf (I)} is satisfied, on the one hand, one has

\[
\liminf \limits_{k \rightarrow + \infty} \mu_k^{\frac{p}{p - q_k}} \geq e^{-\frac np}\; .
\]

\n On the other hand, by (\ref{inf}), one gets

\begin{equation}\label{sub}
B(y_{q_k}, A_{q_k}^{\frac{1}{p}}) \subset B(y_{q_k}, \varepsilon d(x_{q_k},y_{q_k}))
\end{equation}

\n for $k$ large enough. Thus, joining Lemma \ref{CA}, (\ref{3int}) and (\ref{pf}), one arrives at the contradiction

\[
0 < e^{-\frac np} \leq \lim_{k \rightarrow + \infty} \int_{B(y_{q_k},A_{q_k}^{\frac{1}{p}})} u_{q_k}^p\; dv_g = 0 \; .
\]

\n Assume then the assertion {\bf (II)}. In this case, for $k$ large, we set

\[
\begin{array}{l}
\tilde{h}_k(x) = g(\exp_{y_k}(A_k^{\frac{1}{p}} x)) \vspace{0,2cm}\\
\psi_k(x) = u_k(y_k)^{-1} u_k(\exp_{y_k}(A_k^{\frac{1}{p}} x))\; .
\end{array}
\]

\n Thanks to (\ref{l3ep}) and (\ref{nu}), we have

\[
\Delta_{p,h_k} \psi_k + A_k C_k \psi_k^{p - 1} + \frac{1 - \theta_k}{\theta_k} \frac{\nu_k}{\mu_k} \psi_k^{q_k - 1} =
\frac{\nu_k}{\theta_k} \psi_k^{p - 1} \hspace{0,2cm} \mbox{on} \hspace{0,2cm} B(0,3) \; .
\]

\n Rewriting this equation as

\[
\Delta_{p,h_k} \psi_k + A_k C_k \psi_k^{p - 1} + \frac{\nu_k}{\theta_k}\left(\frac{1 - \theta_k}{\mu_k} -1\right)
\psi_k^{q_k - 1} = \frac{\nu_k}{\theta_k}\left( \psi_k^{p - 1} - \psi_k^{q_k - 1}\right) \hspace{0,2cm} \mbox{on} \hspace{0,2cm} B(0,3) \; ,
\]

\n from the mean value theorem, one gets

\begin{equation} \label{norm}
\Delta_{p,h_k} \psi_k + A_k C_k \psi_k^{p - 1} + \frac{\nu_k}{\theta_k}\left(\frac{1 - \theta_k}{\mu_k} -1\right) \psi_k^{q_k - 1} = \frac{\nu_k (pq_k - nq_k + np)}{n}
\psi_k^{\rho_k} \log (\psi_k)  \hspace{0,2cm} \mbox{on} \hspace{0,2cm} B(0,3)
\end{equation}

\n for some $\rho_k \in (q_k-1,p-1)$.

\n Using (\ref{inf}), (\ref{lim}) and the fact that $\frac{1 - \theta_k}{\mu_k} - 1 > 0$ for $k$ large, one easily deduces that $\psi_k \rightarrow \psi$ in
$W^{1,p}(B(0,2))$ and, by a Moser's iteration, one has $\psi \not \equiv 0$.

\n Let $h \in C^1_0(B(0,2))$ be a fixed cutoff function such that $h \equiv 1$ in $B(0,1)$ and $h \geq 0$. Taking $\psi h^p$ as a test function in (\ref{norm}), one obtains

\[
\limsup_{k \rightarrow + \infty} \frac{1}{\theta_k}\left(\frac{1 - \theta_k}{\mu_k} - 1\right) < c\; .
\]

\n Therefore, up to a subsequence, we can write

\[
\lim_{k \rightarrow + \infty} \frac{1}{\theta_k}\left(\frac{1 -
\theta_k}{\mu_k} - 1\right) = \gamma \geq 0\; ,
\]

\n so that $\mu_k \rightarrow 1$. Using the definition of $\theta_k$, we then derive

\[
\lim \limits_{k \rightarrow + \infty} \mu_k^{\frac{p}{p - q_k}} = e^{-(1 + \gamma)\frac{n}{p}}\; .
\]

\n At last, combining Lemma \ref{CA}, (\ref{3int}), (\ref{pf}) and the above limit, we obtain the contradiction

\[
0 <  e^{-(1 + \gamma)\frac{n}{p}} \leq \lim_{k \rightarrow + \infty} \int_{B(y_k,A_k^{\frac{1}{p}})} u_k^p\; dv_g = 0\; .
\]

\bl\\

Finally, we turn our attention to the final argument of the proof of Theorem \ref{lB}. We recall that our goal is to prove that the sequence $(C_k)$ is bounded by assuming that $A_k \rightarrow 0$ as $k \rightarrow + \infty$ (situation {\bf (ii)}). This step consists of several integral estimates around the maximum point $x_k$ of $u_k$ and Lemma \ref{UE} plays a central role on some of them.

In what follows, several possibly different positive constants independent of $k$ will be denoted by $c$ or $c_1$.

Assume, without loss of generality, that the radius of injectivity of $M$ is greater than $2$. Let $\eta \in C^1_0(\R)$ be a cutoff function such that $\eta = 1$ on $[0,1)$, $\eta = 0$ on $[2, \infty)$ and $0 \leq \eta \leq 1$ and define $\eta_k(x) = \eta(d_g(x,x_k))$.

The sharp Euclidean $L^p$-Nash inequality asserts that

\[
\left( \int_{B(0,2)} (u_k \eta_k)^p\; dx \right)^{\frac{1}{\theta_k}} \leq N(p,q_k) \left( \int_{B(0,2)} |\nabla (u_k \eta_k)|^p\; dx\right) \left(\int_{B(0,2)} (u_k \eta_k)^{q_k}\; dx \right)^{\frac{p(1 - \theta_k)}{q_k \theta_k}}\; .
\]

\n Expanding the metric $g$ in normal coordinates around $x_k$, one locally gets

\[
(1 - c d_g(x,x_k)^2) dv_g \leq dx \leq (1 + c d_g(x,x_k)^2) dv_g
\]

\n and

\[
|\nabla(u_k \eta_k)|^p \leq |\nabla_g(u_k \eta_k)|^p (1 + c \; d_g(x,x_k)^2) \; .
\]

\n Thus,

\[
\left( \int_{B(0,2)} u_k^p \eta_k^p\; dx \right)^{\frac{1}{\theta_k}} \leq \left( N(p,q_k) A_k \int_{B(x_k,2)} |\nabla_g (u_k \eta_k)|^p\; dv_g + c A_k \int_{B(x_k,2)} |\nabla_g (u_k\eta_k)|^p d_g(x,x_k)^2\; dv_g \right)
\]

\[
\times \left( \frac{\int_{B(0,2)} u_k^{q_k} \eta_k^{q_k}\; dx}{\int_M u_k^{q_k}\; dv_g} \right)^{\frac{p(1 - \theta_k)}{q_k \theta_k}} \; .
\]

\n Using now the inequality

\[
|\nabla_g (u_k \eta_k)|^p \leq |\nabla_g u_k|^p \eta_k^p + c |\eta_k \nabla_g u_k|^{p - 1} |u_k \nabla_g \eta_k| + c |u_k \nabla_g \eta_k|^p \; ,
\]

\n and denoting

\[
X_k =  A_k \int_M \eta_k^p |\nabla_g u_k|^p d_g(x,x_k)^2\; dv_g
\]

\n and

\[
Y_k = A_k \int_M |\nabla_g u_k|^{p -1} |\nabla_g \eta_k| u_k\; dv_g \; ,
\]

\n we deduce that

\begin{equation}\label{af1}
\left( \int_{B(0,2)} u_k^p \eta_k^p\; dx \right)^{\frac{1}{\theta_k}} \leq \left( N(p,q_k) A_k \int_M |\nabla_g u_k|^p \eta_k^p\; dv_g + c X_k + c Y_k + c A_k \right)
\left( \frac{\int_{B(0,2)} u_k^{q_k} \eta_k^{q_k}\; dx}{\int_M u_k^{q_k}\; dv_g} \right)^{\frac{p(1 - \theta_k)}{q_k \theta_k}}\; .
\end{equation}

\n On the other hand, choosing $u_k \eta_k^p$ as a test function in (\ref{l3ep}) and using (\ref{nu}), one gets

\[
N(p,q_k) A_k \int_M |\nabla_g u_k|^p \eta_k^p\; dv_g \leq 1 - N(p,q_k) A_k C_k + \frac{1}{\theta_k} \left( \int_M u_k^p \eta_k^p\; dv_g - \frac{\int_M u_k^{q_k} \eta_k^p\; dv_g}{\int_M u_k^{q_k}\;
dv_g} \right)
\]

\[
+ c \int_M |\nabla_g u_k|^{p -1} |\nabla_g \eta_k| u_k\; dv_g \; .
\]

\n Using (\ref{ldf}), Proposition \ref{P.1} and Lemma \ref{UE} with a suitable value of $\lambda$, the last integral can be estimated as

\[
\int_M |\nabla_g u_k|^{p -1} |\nabla_g \eta_k| u_k\; dv_g \leq c \left( \int_M |\nabla_g u_k|^p\; dv_g \right)^{\frac{p-1}{p}} \left( \int_{B(x_k,2) \setminus B(x_k,1)} u_k^p\; dv_g \right)^{\frac{1}{p}} \leq c A_k\; ,
\]

\n so that

\begin{equation}\label{af2}
N(p,q_k) A_k \int_M |\nabla_g u_k|^p \eta_k^p\; dv_g \leq 1 - c_1 A_k C_k + \frac{1}{\theta_k} \left( \int_M u_k^p \eta_k^p\; dv_g - \frac{\int_M u_k^{q_k} \eta_k^p\; dv_g}{\int_M u_k^{q_k}\;
dv_g} \right) + c A_k
\end{equation}

\n for $k$ large enough.

\n Let

\[
Z_k = \frac{1}{\theta_k} \left( \int_M u_k^p \eta_k^p\; dv_g - \frac{\int_M u_k^{q_k} \eta_k^{q_k}\; dv_g}{\int_M u_k^{q_k}\; dv_g} \right) \; .
\]

\n By the mean value theorem and Lemma \ref{UE}, there exists $\gamma_k \in (q_k,p)$ such that

\begin{equation}\label{af3}
\left| Z_k - \frac{1}{\theta_k}\left( \int_M u_k^p \eta_k^p\; dv_g - \frac{\int_M u_k^{q_k} \eta_k^p\; dv_g}{\int_M u_k^{q_k}\; dv_g} \right)
\right| \leq \frac{1}{\theta_k} \left| \frac{\int_M u_k^{q_k} (\eta_k^{q_k} - \eta_k^p)\; dv_g}{\int_M u_k^{q_k}\; dv_g} \right| \leq \frac{pq_k - nq_k + np}{n} \frac{\int_M \eta_k^{\gamma_k} | \log \eta_k |
u_k^{q_k}\; dv_g}{\int_M u_k^{q_k}\; dv_g}
\end{equation}

\[
\leq c \frac{\int_{B(x_k,2) \setminus B(x_k,1)} u_k^{q_k}\; dv_g}{\int_M u_k^{q_k}\; dv_g} \leq c A_k
\]

\n for $k$ large enough.

\n Plugging (\ref{af3}) into (\ref{af2}) and after (\ref{af2}) into (\ref{af1}), one arrives at

\begin{equation} \label{af4}
\left( \int_{B(0,2)} u_k^p \eta_k^p\; dx \right)^{\frac{1}{\theta_k}} \leq \left( 1 - c A_k C_k + Z_k + c X_k + c Y_k + c A_k\right) \left( \frac{\int_{B(0,2)} u_k^{q_k} \eta_k^{q_k}\; dx}{\int_M u_k^{q_k}\; dv_g} \right)^{\frac{p(1 - \theta_k)}{q_k \theta_k}}
\end{equation}

\n for $k$ large enough.

\n In order to estimate $X_k$, we take $u_k d_g^2 \eta_k^p$ as a test function in (\ref{l3ep}). From this choice, we derive

\[
X_k \leq \frac{\nu_k}{\theta_k}\left( \int_M u_k^p \eta_k^p d_g(x,x_k)^2\; dv_g - \frac{\int_M u_k^{q_k} \eta_k^{q_k} d_g(x,x_k)^2\; dv_g}{\int_M u_k^{q_k}\;
dv_g}\right) + c A_k \int_M u_k \eta_k^p |\nabla_g u_k|^{p - 1} d_g(x,x_k)\; dv_g
\]

\[
+ c \frac{\int_M u_k^{q_k} \eta_k^{q_k} d_g(x,x_k)^2\; dv_g}{\int_M u_k^{q_k} dv_g} + c Y_k + c A_k\; .
\]

\n We now estimate the first two terms of the right-hand side above. Namely, after a change of variable, one has

\[
\frac{1}{\theta_k} \int_{M} \left| u_k^p \eta_k^p d_g(x,x_k)^2 - \frac{u_k^{q_k} \eta_k^{q_k} d_g(x,x_k)^2}{\int_M u_k^{q_k}\; dv_g}\right|\; dv_g \leq c A_k^{\frac{2}{p}} \frac{1}{\theta_k} \int_{B(0,2 A_k^{- \frac{1}{p}})} \left| \varphi_k^p \tilde{\eta}_k^p - \varphi_k^{q_k} \tilde{\eta}_k^{q_k} \right| |x|^2\; dx
\]

\[
= c A_k^{\frac{2}{p}} \int_{B(0,2 A_k^{- \frac{1}{p}})} (\varphi_k \tilde{\eta}_k)^{\rho_k} \left| \log(\varphi_k \tilde{\eta}_k) \right| |x|^2\; dx
\]

\n for some $\rho_k \in (q_k,p)$, where $\tilde{\eta}_k(x) = \eta_k(A_k^{\frac{1}{p}} x)$.

\n Thus, using Lemma \ref{UE} and the assumption $p \leq 2$, one obtains

\begin{equation} \label{est6}
\frac{1}{\theta_k} \int_{M} \left| u_k^p \eta_k^p d_g(x,x_k)^2 - \frac{u_k^{q_k} \eta_k^{q_k} d_g(x,x_k)^2}{\int_M u_k^{q_k}\; dv_g}\right|\; dv_g \leq c A_k\; .
\end{equation}

\n In particular, since

\[
\int_M u_k^p \eta_k^p d_g(x,x_k)^2\; dv_g - \frac{\int_M u_k^{q_k} \eta_k^{q_k} d_g(x,x_k)^2\; dv_g}{\int_M u_k^{q_k}\; dv_g} = \int_{M} \left( u_k^p \eta_k^p d_g(x,x_k)^2 - \frac{u_k^{q_k} \eta_k^{q_k} d_g(x,x_k)^2}{\int_M u_k^{q_k} \; dv_g}\; \right) dv_g\; ,
\]

\n we have

\[
\frac{1}{\theta_k}\left| \int_M u_k^p \eta_k^p d_g(x,x_k)^2\; dv_g - \frac{\int_M u_k^{q_k} \eta_k^{q_k} d_g(x,x_k)^2\; dv_g}{\int_M u_k^{q_k}\;
dv_g}\right| \leq c A_k
\]

\n for $k$ large enough.

\n Besides, thanks to (\ref{ldf}), Proposition \ref{P.1}, Lemma \ref{UE} and the fact that $p \leq 2$, we derive

\[
\int_M u_k \eta_k^p |\nabla_g u_k|^{p - 1} d_g(x,x_k)\; dv_g \leq c \left( \int_M |\nabla_g u_k|^p\; dv_g \right)^{\frac{p-1}{p}} \left( \int_{B(x_k,2)} u_k^p d_g(x,x_k)^p\; dv_g \right)^{\frac{1}{p}} \]

\[
\leq c A_k^{\frac{2-p}{p}} \left( \int_{B(0,2 A_k^{- \frac{1}{p}})} \varphi_k^p |x|^p\; dx \right)^{\frac{p-1}{p}} \leq c
\]

\n for $k$ large enough.

\n So, the above estimates guarantees that

\[
X_k \leq  c \frac{\int_M u_k^{q_k} \eta_k^{q_k} d_g(x,x_k)^2\; dv_g}{\int_M u_k^{q_k}\; dv_g} + c Y_k + c A_k\; .
\]

\n Evoking again Lemma \ref{UE} and the condition $p \leq 2$, one has

\begin{equation}\label{f1}
\frac{\int_M u_k^{q_k} \eta_k^{q_k} d_g(x,x_k)^2\; dv_g}{\int_M u_k^{q_k}\; dv_g} \leq c A_k^{\frac{2}{p}}  \int_{B(0,2A_k^{-1/p})} \varphi_k^{q_k} |x|^2\; dh_k \leq c A_k
\end{equation}

\n Also, it follows directly from (\ref{ldf}) and Proposition \ref{P.1} that

\begin{equation}\label{est3}
Y_k \leq c A_k^{\frac{1}{p}} \int_{B(x_k,2) \backslash B(x_k,1)} u_k^p\; dv_g \leq c A_k\; ,
\end{equation}

\n so that

\begin{equation}\label{f2}
X_k \leq c A_k\; .
\end{equation}

\n Thus, plugging (\ref{est3}) and (\ref{f2}) into (\ref{af4}), one gets

\[
\left(\int_{B(0,2)} u_k^p \eta_k^p\; dx \right)^{\frac{1}{\theta_k}} \leq \left(1 + Z_k - c_1 A_k C_k + c A_k \right) \left(\frac{\int_{B(0,2)} (u_k \eta_k)^{q_k}\; dx}{\int_M u_k^{q_k}\; dv_g} \right)^{\frac{p(1 - \theta_k)}{q_k \theta_k}}
\]

\n for $k$ large enough.

\n Taking logarithm of both sides and using the fact that $\frac{p(1 - \theta_k)}{q \theta_k} = \frac{1}{\theta_k} - \frac{n - p}{n}$, one has

\begin{equation}\label{d1}
\frac{1}{\theta_k} \left( \log \int_{B(0,2)} u_k^p \eta_k^p dx - \log \left( \frac{\int_{B(0,2)} u_k^{q_k} \eta_k^{q_k} dx}{\int_M u_k^{q_k}\; dv_g}
\right) \right) \leq \log (1 + Z_k - c_1 A_k C_k + c A_k) - \frac{n - p}{n} \log \left( \frac{\int_{B(0,2)} u_k^{q_k} \eta_k^{q_k}\; dx}{\int_M u_k^{q_k}\; dv_g}
\right)\; .
\end{equation}

\n By the mean value theorem,

\begin{equation}\label{d2}
\log \int_{B(0,2)} u_k^p \eta_k^p\; dx - \log \left( \frac{\int_{B(0,2)} u_k^{q_k} \eta_k^{q_k}\; dx}{\int_M u_k^{q_k}\; dv_g} \right)
= \frac{1}{\tau_k} \left( \int_{B(0,2)} u_k^p \eta_k^p\; dx - \frac{\int_{B(0,2)} u_k^{q_k} \eta_k^{q_k}\; dx}{\int_M u_k^{q_k}\; dv_g} \right)
\end{equation}

\n for some number $\tau_k$ between the expressions

\[
\int_{B(0,2)} u_k^p \eta_k^p\; dx\ \ {\rm and}\ \ \frac{\int_{B(0,2)} u_k^{q_k} \eta_k^{q_k}\; dx}{\int_M u_k^{q_k}\; dv_g} \; .
\]

\n Using Cartan's expansion of $g$ in normal coordinates around $x_k$ and Lemma \ref{UE}, one obtains

\begin{equation} \label{est4}
\max\left\{ \left| \int_{B(0,2)} u_k^p \eta_k^p\; dx - \int_M u_k^p \eta_k^p\; dv_g \right|, \left| \frac{\int_{B(0,2)} u_k^{q_k} \eta_k^{q_k}\; dx}{\int_M u_k^{q_k}\; dv_g} - \frac{\int_M u_k^{q_k} \eta_k^{q_k} dv_g}{\int_M u_k^{q_k} dv_g} \right|\right\} \leq c A_k
\end{equation}

\n for $k$ large enough. Indeed, since $p \leq 2$,

\[
\left| \int_{B(0,2)} u_k^p \eta_k^p\; dx - \int_M u_k^p \eta_k^p\; dv_g \right| \leq c \int_M u_k^p \eta_k^p d_g(x,x_k)^2 \; dv_g \leq c A_k^{\frac{2}{p}}  \int_{B(0,2A_k^{-1/p})} \varphi_k^p |x|^2\; dv_{h_k} \leq c A_k
\]

\n and, by (\ref{f1}),

\[
\left| \frac{\int_{B(0,2)} u_k^{q_k} \eta_k^{q_k}\; dx}{\int_M u_k^{q_k}\; dv_g} - \frac{\int_M u_k^{q_k} \eta_k^{q_k} dv_g}{\int_M u_k^{q_k} dv_g} \right| \leq c \frac{\int_M u_k^{q_k} \eta_k^{q_k} d_g(x,x_k)^2\; dv_g}{\int_M u_k^{q_k}\; dv_g} \leq c A_k\; .
\]

\n Moreover, we also have

\begin{equation} \label{est5}
\max\left\{ \left| \int_M (u_k \eta_k)^p\; dv_g - 1\right|, \left|\frac{\int_M u_k^{q_k} \eta_k^{q_k}\; dv_g}{\int_M u_k^{q_k}\; dv_g} - 1 \right|\right\} \leq c A_k
\end{equation}

\n for $k$ large enough. In fact, by Lemma \ref{UE},

\[
\left| \int_M u_k^p \eta_k^p\; dv_g - 1\right| = \left| \int_M u_k^p \eta_k^p\; dv_g - \int_M u_k^p \; dv_g\right| \leq c \int_{M \setminus B(x_k,1)} u_k^p \; dv_g \leq c A_k
\]

\n and

\[
\left|\frac{\int_M u_k^{q_k} \eta_k^{q_k}\; dv_g}{\int_M u_k^{q_k}\; dv_g} - 1 \right| \leq c \frac{\int_{M \setminus B(x_k,1)} u_k^{q_k}\; dv_g}{\int_M u_k^{q_k}\; dv_g} \leq c A_k\; .
\]

\n Thanks to (\ref{est4}) and (\ref{est5}), one easily deduces that $\tau_k^{-1} = 1 + O(A_k)$. Then, by (\ref{d2}),

\begin{equation}\label{d3}
\frac{1}{\theta_k} \left( \log \int_{B(0,2)} u_k^p \eta_k^p\; dx - \log \left( \frac{\int_{B(0,2)} u_k^{q_k} \eta_k^{q_k}\; dx}{\int_M u_k^{q_k}\; dv_g} \right) \right)
= \frac{1}{\theta_k} \left( \int_{B(0,2)} u_k^p \eta_k^p\; dx - \frac{\int_{B(0,2)} u_k^{q_k} \eta_k^{q_k}\; dx}{\int_M u_k^{q_k}\; dv_g} \right) (1 + O(A_k))\; .
\end{equation}

\n But, by Cartan's expansion and (\ref{est6}), we have

\[
\frac{1}{\theta_k} \left( \int_{B(0,2)} u_k^p \eta_k^p\; dx - \frac{\int_{B(0,2)} u_k^{q_k} \eta_k^{q_k}\; dx}{\int_M u_k^{q_k}\; dv_g} \right) = \frac{1}{\theta_k} \left(\int_{B(0,2)} u_k^p \eta_k^p - \frac{u_k^{q_k} \eta_k^{q_k}}{\int_M u_k^{q_k}\; dv_g}\; dx \right)
\]

\[
=  \frac{1}{\theta_k} \left(\int_M u_k^p \eta_k^p - \frac{u_k^{q_k} \eta_k^{q_k}}{\int_M u_k^{q_k}\; dv_g}\; dv_g \right) + \frac{1}{\theta_k} \left(\int_M \left( u_k^p \eta_k^p - \frac{u_k^{q_k} \eta_k^{q_k}}{\int_M u_k^{q_k}\; dv_g} \right) O(d_g(x, x_k)^2)\; dv_g \right)
\]

\[
= Z_k + O(A_k)
\]

\n for $k$ large enough.

\n Replacing this inequality in (\ref{d3}), one obtains

\[
\frac{1}{\theta_k} \left(\log \int_{B(0,2)} u_k^p \eta_k^p\; dx - \log \left( \frac{\int_{B(0,2)} u_k^{q_k} \eta_k^{q_k}\; dx}{\int_M u_k^{q_k}\; dv_g} \right)\right) \geq Z_k - c A_k\; .
\]

\n In turn, plugging the above inequality in (\ref{d1}), one has

\[
Z_k - c A_k \leq \log (1 + Z_k - c_1 A_k C_k + c A_k) + \frac{n - p}{n} \left| \log \left( \frac{\int_{B(0,2)} u_k^{q_k} \eta_k^{q_k}\; dx}{\int_M u_k^{q_k}\; dv_g} \right) \right|\; .
\]

\n Finally, using Cartan's expansion of $g$ in normal coordinates, Taylor's expansion of the function $\log$ and Lemma \ref{UE}, one gets

\[
\left| \log \left( \frac{\int_{B(0,2)} u_k^{q_k} \eta_k^{q_k}\; dx}{\int_M u_k^{q_k}\; dv_g} \right) \right| \leq c \frac{\int_{M \backslash B(x_k,1)} u_k^{q_k}\; dv_g}{\int_M u_k^{q_k}\; dv_g} + c \frac{\int_M u_k^{q_k} \eta_k^{q_k} d_g(x,x_k)^2\; dv_g}{\int_M u_k^{q_k}\; dv_g} \leq c A_k\; .
\]

\n In short, for a certain constant $c > 0$, we deduce that

\[
Z_k \leq \log (1 + Z_k - c_1 A_k C_k + c A_k) + c A_k
\]

\n for $k$ large enough.

\n Since $\log x \leq x - 1$ for all $x > 0$, we find

\[
Z_k \leq Z_k - c_1 A_k C_k + c A_k
\]

\n for $k$ large enough, so that the sequence $(C_k)$ is bounded and the proof of Theorem \ref{lB} follows.\bl \\

{\bf Acknowledgements.} The authors are indebted to the referee for his valuable suggestions and comments pointed out concerning this work. The first author was partially supported by CAPES through INCTmat and the second one was partially supported by CNPq and Fapemig.

\end{document}